\documentclass[conference]{IEEEtran}
\IEEEoverridecommandlockouts
\usepackage{cite}
\usepackage{amsmath,amssymb,amsfonts}
\usepackage{algorithmic}
\usepackage{graphicx}
\usepackage{textcomp}
\usepackage{xcolor}
\def\BibTeX{{\rm B\kern-.05em{\sc i\kern-.025em b}\kern-.08em
    T\kern-.1667em\lower.7ex\hbox{E}\kern-.125emX}}
\begin{document}

\title{On extending the class of convex functions\\
}

\author{\IEEEauthorblockN{Shravan Mohan}
\thanks{This work is not related to the work of the author at Amazon.}
\IEEEauthorblockA{Mantri Residency, Bengaluru\\
shravan8587@gmail.com}
}

\maketitle

\begin{abstract}
In this brief note, it is shown that the function $p^\top W \log(p)$ is convex in $p$ if $W$ is a diagonally dominant positive definite M-matrix. The techniques used to prove convexity are well-known in  linear algebra and essentially involves factoring the Hessian in a way that is amenable to martix analysis. Using  similar techniques, two classes of convex homogeneous polynomials is derived - namely, $p^\top W p^2$ and $(p^k)^\top W p^k$ - the latter also happen to be SOS-convex. Lastly, usign the same techniques, it is also shown that the function  $p^\top W e^p$ is convex over the positive reals only if $W$ is a non-negative diagonal matrix. Discussions regarding the utility of these functions and examples accompany the results presented.
\end{abstract}

\begin{IEEEkeywords}
Convexity, entropy, posynomials, Schur complement, SOS-convexity
\end{IEEEkeywords}

\section{Introduction}
\noindent Convex functions play a pivotal role in optimization theory due to their unique properties that guarantee the local minimum of a convex function is also a global minimum. This attribute simplifies optimization problems significantly, making convex functions integral in fields such as operations research, economics, and machine learning, where finding optimal solutions is crucial. Moreover, the structural characteristics of convex functions, such as their ability to approximate complex functions using convex envelopes, contribute to their widespread application in solving large-scale optimization problems efficiently \cite{boyd2004convex, grant2014cvx}. \\\\
The aim of this work is to enrich the set of well-known convex functions with a few more. To the best of the authors knowledge, these functions have not been shown to be convex before. The first in these relates to the entropy function. The entropy function commonly represented in the context of information theory as 
\begin{align}
H(p) = -\sum_{x} p(x) \log p(x), 
\end{align}  
where \( p(x) \) denotes the probability distribution over discrete random variables, serves as a fundamental measure of uncertainty or randomness in a system. This function is notably convex, which can be demonstrated through its second derivative. The strict convexity of the entropy function plays a pivotal role in determining the characteristics of probability distributions that maximize entropy. Specifically, it ensures that there is a unique distribution that achieves this maximum — namely, the uniform distribution. This uniqueness arises because strict convexity guarantees that any local extremum of the entropy function is also the global maximum, which, in the case of entropy, corresponds to the scenario where outcomes are equiprobable, as seen in the uniform distribution. Furthermore, the convexity of the entropy function directly influences the properties of the Kullback-Leibler divergence, which is a measure of how one probability distribution diverges from another. The convexity of KL divergence with respect to its arguments ensures that the divergence is minimized precisely when the two distributions being compared are identical \cite{ash2012information}. \\

\noindent In the first part, the convexity of a novel entropy-like function: \( p^\top W \log(p) \) is established for $W$ belonging to a specific class of matrices. This work also generalizes the standard entropy function, which emerges as a special case when \( W \) is set to the identity matrix. Although the immediate applications of these results are not yet fully clear, the establishment of such properties invites further exploration into their potential utility across various domains. The method employed to prove the convexity of the aforementioned entropy-like function involves factoring the Hessian and making it amenable to well-known techniques from linear algebra, so that sufficient conditions for its positive semi-definiteness can be obtained.\\ 

\noindent The techniques used are also extendable to a range of polynomial functions (specifically, posynomials \cite{boyd2007tutorial}), thereby further enriching the set of known convex functions. This set of results forms the second part of this paper. Incidentally, the class of polynomials which are shown to be convex also happen to be SOS-convex, a well-known subset of all convex polynomials \cite{ahmadi2013complete}. The last result in this paper shows that the function $x^\top W e^x$ is convex over the positive reals only if $W$ is a positive diagonal matrix. The function $xe^x$ appears in several applied areas, such as, in solutions to ODEs \cite{kailath1980linear}, the Lambert W function \cite{mezo2022lambert}. In other words, this result negates the possibility of having mixed terms ($x_ie^{x_j},~i\neq j$) with the additional requirement of convexity.  \\

\noindent The next section presents the main results of this paper along with their detailed proofs, demonstrating the application of the aforementioned linear algebraic techniques to establish convexity. Each theorem is supplemented with a discussion for clarity. This is followed by concluding remarks that summarize the key findings and suggest directions for future research.\\

\section{Preliminaries}
\noindent The term $\mathcal{R}$ represents the set of real numbers and $\mathcal{R}^+$ represents the set of positive reals. Matrices and vectors are  denoted by uppercase  and lowercase  letters, respectively. The entries of a matrix \(A\) are represented as \(A_{i,j}\), where \(i\) and \(j\) index rows and columns. Matrix operations used here include the transpose \(A^\top\), inverse \(A^{-1}\), and the Hadamard product \(A \circ B\). The term $D_v$ refers to a diagonal matrix with $v$ as its diagonal elements. The inequality $A \geq 0$ implies that all elements of $A$ are non-negative, while the inequality $A\succeq 0$ implies that the matrix $A$ is positive semidefinite.\\

\noindent \textbf{\underline{Convex function \cite{boyd2004convex}}}: A twice differentiable function \( f: \mathbb{R}^n \to \mathbb{R} \) is defined to be convex on a domain \( D \subseteq \mathbb{R}^n \) if the Hessian matrix of \( f \), denoted as \( \nabla^2 f(x) \), is positive semidefinite for all \( x \) in \( D \). \\

\noindent \textbf{\underline{M-matrix \cite{poole1974survey}}}: A matrix \( A \in \mathbb{R}^{n \times n} \) is called an M-matrix if it can be expressed in the form 
\begin{align}
A = sI - B,
\end{align}
where \( I \) is the identity matrix, \( B \) is a nonnegative matrix (i.e., all entries of \( B \) are nonnegative), and \( s \) is a scalar greater than or equal to the spectral radius of \( B \). \\

\noindent \textbf{\underline{Inverse-Positivity of M-Matrices \cite{poole1974survey}}}: Let \( A \) be an M-matrix. Then the inverse of \( A \) is nonnegative, i.e., 
\begin{align}
A^{-1} \geq 0,
\end{align}
where the inequality is element-wise.\\

\noindent\textbf{ \underline{Sylvester's Law of Inertia \cite{bhatia2009positive}}}: 
Let \( A \) be a real symmetric matrix. Then, for any invertible matrix \( P \), the matrix \( B = P^\top A P \) is congruent to \( A \) and has the same inertia. That is, the number of positive eigenvalues, the number of negative eigenvalues, and the number of zero eigenvalues of \( A \) remain invariant under such transformations.\\

\noindent \textbf{\underline{Schur's Complement Theorem \cite{bhatia2009positive}}}: Let \( M \) be a block matrix given by
\begin{align}
M = \begin{bmatrix}
A & B \\
C & D
\end{bmatrix}
\end{align}
where \( A \), \( B \), \( C \), and \( D \) are matrices of appropriate dimensions and \( A \) is invertible. The Schur complement of \( A \) in \( M \) is defined as
\begin{align}
M/A = D - CA^{-1}B.
\end{align}
The matrix \( M \) is positive definite if and only if both \( A \) and the Schur complement \( M/A \) are positive definite.\\

\noindent \textbf{\underline{Schur Product Theorem:}} Let $A$ and $B$ be two $n \times n$ complex Hermitian positive semidefinite matrices. The Hadamard product $(A \circ B)$ is also positive semidefinite.\\

\noindent \textbf{\underline{Posynomial \cite{boyd2007tutorial} }}: A polynomial defined over the set of positive reals is called a posynomial if all coefficients are non-negative. \\

\noindent \textbf{\underline{Sum of Squares(SOS)-convexity \cite{ahmadi2013complete}}}: A polynomial $f(x)$ is SOS-convex if and only if $y^\top \nabla^2(f) y$ is a sum of squares polynomial in $(x, y)$. 
\section{The Main Results}

\noindent \textbf{\underline{Theorem 1}}: The function $p^\top W \log(p)$ is convex in $p \in \mathcal{R}_+^n$  if the matrix $W$ is a diagonally dominant positive definite M-matrix. Here, $W\in \mathcal{R}^{n\times n}$ and $\log(.)$ is applied element-wise.\\

\noindent \textbf{\underline{Proof}}:
Suppose $f(p) = g^\top(p)h(p)$, where $g$ and $h$ are smooth functions. Then,
\begin{align}
\frac{\partial f}{\partial p_i} = \frac{\partial g}{\partial p_i}^\top h + g^\top \frac{\partial h}{\partial p_i}.
\end{align}
And hence, 
\begin{equation}\small
\frac{\partial^2 f}{\partial p_i \partial p_j} = \frac{\partial^2 g}{\partial p_i \partial p_j}^\top h + \frac{\partial g}{\partial p_i}^\top\frac{\partial h}{\partial p_j} + \frac{\partial g}{\partial p_j}^\top\frac{\partial h}{\partial p_i} + g^\top \frac{\partial^2 h}{\partial p_i \partial p_j}.
\end{equation}
Now, substituting 
\begin{align}
g = Wp ~~\&~~ h = \log(p),
\end{align}
we have that 
\begin{equation}
\begin{aligned}
    \frac{\partial^2 g}{\partial p_i \partial p_j}^\top h &~~=~~ 0, ~~\forall i,j,\\
    \frac{\partial g}{\partial p_i}^\top\frac{\partial h}{\partial p_j} &~~=~~ \frac{e_i^\top W e_j^\top}{p_j},  ~~\forall i,j,\\   \frac{\partial g}{\partial p_j}^\top\frac{\partial h}{\partial p_i} &~~=~~ \frac{e_j^\top W e_i^\top}{p_i},  ~~\forall i,j,\\
    g^\top \frac{\partial^2 h}{\partial p_i \partial p_j} &~~=~~ \begin{cases} 
    0, & \text{if } i\neq j, \\
    -\frac{p^\top W e_i}{p_i^2}, & \text{if } i=j.
    \end{cases}
\end{aligned}
\end{equation}
Simplifying further, we get that $\nabla^2(f)$ is equal to
\begin{equation}
\begin{aligned}
     &WD_{\frac{1}{p}} + D_{\frac{1}{p}} W - D_{\frac{1}{p}} D_{Wp} D_{\frac{1}{p}}\\
    =&D_{\frac{1}{p}} \left(D_{p}W + WD_{p} - D_{Wp}\right) D_{\frac{1}{p}}.
\end{aligned} 
\end{equation}
By Sylvester's Intertia Theorem, it is clear that 
\begin{align}
    \nabla^2(f) \succeq 0 \rightleftharpoons \left(D_{p}W + WD_{p} - D_{Wp}\right) \succeq 0.
\end{align}
For convexity, the LMI above must be true for all $p$. This implies that it must be true for $p=e_i,~\forall i$. Thus,
\begin{align}
    \left(D_{e_i}W + WD_{e_i} - D_{We_i}\right) \succeq 0.
\end{align}
To break this down further, consider the case $i = 1$. Let the first row of $W$ be given by $[u~~v]$, where $u$ is the first element of the row and $v$ is the vector representing the rest of the elements. In that case, the LMI can be written as 
\begin{align}
    \begin{bmatrix}
        u & v \\
        v^\top & -D_v
    \end{bmatrix} \succeq 0.
\end{align}
By Schur's complement, the LMI is true if and only if
\begin{equation}
\begin{aligned}
    v &\leq 0, \text{~~and} \\
    u + \sum_{i=1}^{n-1} v &\geq 0.
\end{aligned}
\end{equation}
The inequalities imply that  the first row of $W$ must be diagonally dominant and its off diagonal elements must be negative. Similar inequalities would also be true for the rest of the rows. Now, if the above inequalities hold for each row, or equivalently $W$ is an M-matrix, $\nabla^2(f)$ would be positive semi-definite at every $p$, implying convexity. $\hfill{\blacksquare}$\\

\noindent \textbf{\underline{Corollary 2}}: The function $p^\top \log(Wp)$ is convex in $p\in \mathcal{R}^n_+$ if $W$ is the inverse of a positive definite diagonally dominant M-matrix. \\

\noindent \textbf{\underline{Proof}}: Consider the linear transformation of variables given by
\begin{align}
q = Wp.
\end{align}
With that, the function $p^\top W \log(p)$ can be written as
\begin{align}
q^\top W^{-1}\log(q).
\end{align}
By the following relation
\begin{align}
\nabla^2_q(f) = W \nabla^2_p(f)|_{W^{-1}q} W,
\end{align}
and using the previous theorem, it is clear that that the function $p\top \log(Wp)$ is convex in $p\in \mathcal{R}^+$ if $W$ is the inverse of a positive definite diagonally dominant M-matrix. $\hfill{\blacksquare}$\\

\noindent The theorem on the convexity of $p^\top W \log(p)$ subsumes the case of standard entropy when $W$ is equal to the identity matrix. However, for a general $W$, this cannot be deemed as entropy in the information theoretic sense, since the function would not be symmetric in $p$; a key requirement for entropy \cite{ash2012information}. Same is the case with the corollary; that too boils down to the standard entropy when $W$ is equal to the identity matrix. Also, since the inverse of a M-matrix is always non-negative, this implies that the function $p^\top \log(Wp)$ is well-defined for $p\in \mathcal{R}^n$. \\

\noindent Next, as mentioned earlier, we use the similar techniques to mine for classes of convex homogeneous polynomials. As will be shown later, these happen to be a subset of posynomials. Posynomials form the backbone of the field of geometric programming. This is a  technique used for solving  optimization problems characterized by the objective and constraint functions being posynomials. This method has been used  in diverse fields such as network design, structural engineering, and resource allocation \cite{boyd2007tutorial}. The optimization problem is not itself convex, but becomes convex via a logarithmic change of variables and the application of logarithmic transformations to the objective and constraints. Take, for example, the function \( f(x, y) = x^2 y^2 \). To examine its convexity, consider the Hessian matrix of \( f \), given by
\begin{align}
H_f = \begin{bmatrix}
\frac{\partial^2 f}{\partial x^2} & \frac{\partial^2 f}{\partial x \partial y} \\
\frac{\partial^2 f}{\partial y \partial x} & \frac{\partial^2 f}{\partial y^2}
\end{bmatrix}
= \begin{bmatrix}
2y^2 & 4xy \\
4xy & 2x^2
\end{bmatrix}.
\end{align}
The determinant of the Hessian matrix, \( \text{det}(H_f) = 4x^2 y^2 - 16x^2 y^2 = -12x^2 y^2 \), is clearly negative unless one of the variables \( x \) or \( y \) is zero, indicating the presence of saddle points and confirming the function's non-convexity in its original form. This analysis highlights a  limitation of geometric programming, which excludes, for example, scenarios with convex quadratic constraints possessing non-negative coefficients. Establishing the convexity of posynomials in their original variables can be useful, as it expands the scope to a broader range of convex optimization problems, beyond those amenable to standard posynomial transformations. \\

\noindent Deciding convexity of polynomials is not an easy problem. As shown in \cite{ahmadi2013np}, deciding convexity of quartics, and of larger degree polynomials is a NP-hard problem. SOS-convexity \cite{ahmadi2013complete}, which is a special case of convexity. In fact, the authors of  \cite{ahmadi2013complete} established that SOS-convexity implies convexity if and only if one of the following conditions is met: (a) \( n = 2 \), (b) \( 2d = 2 \), or (c) \( n = 3 \) and \( 2d = 4 \). Importantly, deciding whether a polynomial is SOS convex can be done in polynomials time and this essentially relies on the fact that deciding the SOS nature of polynomials is decidable in polynomial time.  The polynomials examined in the first result does not satisfy this criteria (the degree of that class polynomials is 3). The second class of polynomials does fall within the bracket of SOS-convex polynomials, but has the additional  advantage of explicit and simple construct.\\

\noindent \textbf{\underline{Theorem 3}}: The function $p^\top W p^2$ is convex in $p\in \mathcal{R}^n_+$ if $W$ is a non-negative symmetric matrix and 
\begin{align}
W_{i,i} \geq \frac{1}{3}\sum_{j\neq i} W_{i,j}, ~\forall ~i. 
\end{align}\\
\noindent \textbf{\underline{Proof}}: The proof method here is same as that of Theorem 1. Again, suppose $f(p) = g^\top(p)h(p)$, where $g$ and $h$ are smooth functions. With $g(p) = Wp$ and $h(p)=p^2$, we have that
\begin{align}
\nabla^2(f) = WD(p) + D(p)W + D(Wp).
\end{align}
First of all, since $W$ is a non-negative matrix, the third term of the $\nabla^2(f)$ is positive semidefinite for all $p\in \mathcal{R}_+^n$. Next, note that the expression for $\nabla^2(f)$ above is positive semidefinite for all $p\in \mathcal{R}_+^n$ if it is positive semidefinite for $p=e_i, ~\forall ~i$. Consider $p=e_1$ and as before, let the first row of $W$ be given by $[u~~v]$, where $u$ is the first element of the row and $v$ is the vector representing the rest of the elements. Then, one obtains the following:  
\begin{align}
   \nabla^2(f)\big|_{e_1} =  \begin{bmatrix}
        3u & v \\
        v^\top & D(v)
    \end{bmatrix}.
\end{align}
Now, by Schur complement, the Hessian at $e_1$ is positive definite if
\begin{align}
u \geq \frac{1}{3}\sum_{i=1}^{n-1}v_i. 
\end{align}
Extending the same analysis for all $e_i$'s, one obtains that under the hypothesis on $W$, the function is indeed convex. $\hfill{\blacksquare}$\\

\noindent Although the above result involves homogeneous polynomials, a set of non-homogeneous convex degree 3 polynomials can be constructed by composing  the function with non-negative affine transformation (using the result that $f(Ax+b)$ is convex if $f(x)$ is convex). Further a degree 3 polynomial cannot be SOS-convex since the Hessian is linear in $p$.  \\

\noindent \textbf{\underline{Theorem 4}}: The function $(p^2)^\top W p^2$ is convex in $p\in \mathcal{R}^n_+$ if $W$ is a non-negative matrix and 
\begin{align}
\Tilde{W} = \begin{bmatrix}
        \frac{3}{2} W_{1,1} & W_{1,2} & \cdots  & W_{1,n}\\
         W_{2,1} & \frac{3}{2} W_{2,2} & \cdots & W_{2,n}\\
         \vdots & \vdots & \vdots &\vdots \\
         W_{n,1} & W{n,2} &\cdots & \frac{3}{2} W_{n,n}
\end{bmatrix} \succeq 0. 
\end{align}
In fact, the function is also SOS-convex.
\\

\noindent \textbf{\underline{Proof}}: As before, suppose $f(p) = g^\top(p)h(p)$, where $g(p) = Wp^2$ and $h(p)=p^2$. This implies that
\begin{align}
\nabla^2(f) = 8 W\circ (p . p^\top) + 4D_{Wp^2}.
\end{align}
Note that the above expression can be recombined as 
\begin{align}
\nabla^2(f) = 12 W\circ (p . p^\top) + 4D_{W_0p^2},
\end{align}
where 
\begin{equation}
\begin{aligned}
W_0 = \begin{bmatrix}
        0 & W_{1,2} & \cdots  & W_{1,n}\\
         W_{2,1} & 0 & \cdots & W_{2,n}\\
         \vdots & \vdots & \vdots &\vdots \\
         W_{n,1} & W{n,2} &\cdots & 0
\end{bmatrix}.
\end{aligned}
\end{equation}
This can further be  rewritten as 
\begin{align}
\nabla^2(f) = 8 \Tilde{W}\circ (p . p^\top) + 4D_{W_0p^2}.
\end{align}
Firstly, since $W_0$ is nonnegative, the second term is positive semidefinite for all $p\in \mathcal{R}^n$. Secondly, note that the second term is the Hadamard product between two positive semidefinite matrices, that is, $\Tilde{W}$ and $p.p^\top$. And since the Hadamard product of two positive semidefinite matrices is also positive semidefinite (by Schur product theorem), it too is positive semidefinite. This implies that the $\nabla^2(f)$ is positive definite for all $p\in \mathcal{R}^n$ and hence the function is convex everywhere.\\

\noindent To show that the function is also SOS-convex, note the following:
\begin{equation}
\begin{aligned}
y^\top \nabla^2(f) y &= 
8 y^\top  \left(\Tilde{W}\circ (p . p^\top)\right) y + 4y^\top D_{W_0p^2} y\\
&=8\text{Tr}\left(D_y\Tilde{W}D_y p p^\top\right) + 4y^\top D_{W_0p^2} y\\
&=8 p^\top D_y\Tilde{W}D_y p  + 4y^\top D_{W_0p^2} y.
\end{aligned}
\end{equation}
It is clear that the above is a SOS polynomial (since $\Tilde{W}$ and $D_{W_0p^2}$ are both positive semidefinite), and hence the function is SOS-convex. 
$\hfill{\blacksquare}$\\

\noindent The analysis used in the proof of Theorem 4 can be used to extend the idea to higher order polynomials and this is what the next theorem shows. However, as shall be shown the domain of convexity changes based on the parity of the degree.\\

\noindent \textbf{\underline{Theorem 5}}: For an even $k$, the function $(p^k)^\top W p^k$ is convex in $p\in \mathcal{R}^n$ if $W$ is a non-negative and 
\begin{align}
\Tilde{W} = \begin{bmatrix}
        \lambda(k) W_{1,1} & W_{1,2} & \cdots  & W_{1,n}\\
         W_{2,1} & \lambda(k) W_{2,2} & \cdots & W_{2,n}\\
         \vdots & \vdots & \vdots &\vdots \\
         W_{n,1} & W{n,2} &\cdots & \lambda(k) W_{n,n}
\end{bmatrix} \succeq 0,
\end{align}
where 
\begin{align}
\lambda(k) = \frac{k^2 + k(k-1)}{k^2}
\end{align}
In fact, the function is also SOS-convex. For odd $k$, the function above is convex and SOS-convex over $R^n_+$.
\\

\noindent \textbf{\underline{Proof}}: As before, suppose $f(p) = g^\top(p)h(p)$, where $g(p) = Wp^k$ and $h(p)=p^k$, we have that
\begin{align}\resizebox{3in}{!}{
$
\nabla^2(f) = 2k^2 W\circ (p^{k-1} . (p^{k-1})^\top) + 2k(k-1)D_{p^{k-2}}D_{Wp^k}.
$
}
\end{align}
Note that the above expression can be recombined as 
\begin{align}\resizebox{3in}{!}{
$
\nabla^2(f) = 2k^2 \Tilde{W}\circ (p^{k-1} . (p^{k-1})^\top) + 2k(k-1)D_{W_0p^{k-2}},
$
}
\end{align}
where 
\begin{align}
W_0 = \begin{bmatrix}
        0 & W_{1,2} & \cdots  & W_{1,n}\\
         W_{2,1} & 0 & \cdots & W_{2,n}\\
         \vdots & \vdots & \vdots &\vdots \\
         W_{n,1} & W{n,2} &\cdots & 0
\end{bmatrix}
\end{align}
Consider the case when $k$ is even. Since $W_0$ is nonnegative, the second term is positive semidefinite over $\mathcal{R}^n$. Since the first term is the Hadamard product between two positive semidefinite matrices ($\Tilde{W}$ and $p.p^\top$), it too is positive semidefinite. This implies that the $\nabla^2(f)$ is positive definite for all $p\in \mathcal{R}^n$.  \\

\noindent For odd $k$, the second term in the expression for $\nabla^2(f)$ is positive semidefinite if $p\in \mathcal{R}^n_+$. And hence, for odd $k$, the function is convex over $\mathcal{R}^n_+$. Finally, using techniques similar to those used in the proof of Theorem 4, it can be shown that the function is also SOS-convex (for both even and odd $k$, on their respective domains). $\hfill{\blacksquare}$\\

\noindent The next result to be discussed is fundamentally a negative one with respect to the function $p^\top W e^p$. As mentioned earlier, this function appears in several fields of mathematics, and therefore a result such as this would prevent the misapplication of convex optimization routines to problems involving this function. This function can also be seen as the composition of the function $p^\top W \log(p)$ with $e^p$. \\

\noindent \textbf{\underline{Theorem 6}}: If $W$ is a symmetric non-negative matrix, the function $p^\top W e^p$ is convex over $\mathcal{R}_+^n$ if and only if $W$ is a diagonal matrix.\\

\noindent \textbf{\underline{Proof}}: As before, suppose $f(p) = g^\top(p)h(p)$, where $g(p) = Wp$ and $h(p)=e^p$. Then, the Hessian of this function boils down to
\begin{equation}
\begin{aligned} 
&\nabla^2(f) = D_{e^p}W + WD_{e^p} + D_{Wp}D_{e^p}\\ 
=& D_{e^{\frac{p}{2}}} \left( D_{e^{\frac{p}{2}}} W D_{e^{-\frac{p}{2}}} + D_{e^{-\frac{p}{2}}} W D_{e^{\frac{p}{2}}} + D_{Wp}\right) D_{e^{\frac{p}{2}}}.
\end{aligned}
\label{eqn:hessian_pWexpp}
\end{equation}
Consider the case where $W$ is a non-diagonal matrix. And without loss of generality, let the first row of $W$ have a non-zero element at column $j$. Consider the point $p = [2\log(M), 0, \cdots, 0]$ for some $M$. Note that principal minor $(1,j)$ of the term within brackets in the expression for Hessian \eqref{eqn:hessian_pWexpp} can be written as:
\begin{align}
\begin{bmatrix}
W_{1,1}(2 + 2\log(M)) & \left(M+\frac{1}{M}\right) W_{1,j} \\
\left(M+\frac{1}{M}\right) W_{1,j} & W_{j,j} + 2W_{1,j}\log(M)
\end{bmatrix}.
\end{align}
It is immediately clear that the determinant of this 2x2 matrix, which is given by
\begin{align} \resizebox{3in}{!}{
$
W_{1,1}(2 + 2\log(M))(W_{j,j} + 2W_{1,j}\log(M)) - \left(M+\frac{1}{M}\right)^2W^2_{1,j}
$
},
\end{align}
becomes negative for a large enough $M$, there by disproving convexity over the entire domain of $\mathcal{R}_+^n$. It is also clear that the only way the function can be convex is when $W$ is a non-negative diagonal matrix. In this case, the function is given by
\begin{align}
\sum_{i=0}^np_ie^{p_i}.
\end{align}
And since the function $pe^p$ is convex over $\mathcal{R}_+^n$, a positive sum of such functions is also convex. 
$\hfill{\blacksquare}$

\section{Conclusion}
\noindent Establishing specific conditions for $f(p)$ under which the function \( p^\top W f(p) \) is convex would be an interesting direction for further research. A general result of this sort would subsume most of the results presented in this paper and would also greatly enhance the known classes of convex functions. Additionally, the conditions for convexity presented in this paper are only sufficient (for example, symmetricity of the matrix $W$). It would be interesting to see if necessary and sufficient conditions can be determined for convexity.
\bibliographystyle{IEEEtran}
\bibliography{biblio}
\end{document}